\documentclass[preprint, 12pt]{elsarticle}
\usepackage{amssymb}
\usepackage{amsmath}
\usepackage{kpfonts}
\usepackage[T1]{fontenc}
\usepackage{xcolor}
\usepackage{hyperref}
\usepackage{mathtools}
\usepackage[utf8]{inputenc}
\usepackage{float} 
\usepackage{natbib}
\usepackage{cleveref}
\newtheorem{con}{Conjecture}
\newtheorem{que}{Question}
\newdefinition{ex}{Example}
\newdefinition{rem}{Remark}
\newdefinition{cons}{Construction}

\begin{document}
\bibliographystyle{abbrv}

\makeatletter
\def\ps@pprintTitle{
\let\@oddhead\@empty
\let\@evenhead\@empty
\let\@evenfoot\@oddfoot}
\makeatother

\begin{frontmatter}

\title{\textbf{On Huisman's conjectures about unramified real curves}}

\author[mario]{MARIO KUMMER}
\ead{mario.kummer@tu-dresden.de}

\author[dimitri]{DIMITRI MANEVICH}
\ead{dimitri.manevich@math.tu-dortmund.de}

\address[mario]{Fakultät Mathematik, Institut für Geometrie, Technische Universität Dresden, Zellescher Weg 12-14, 01062 Dresden, Germany}
\address[dimitri]{Fakultät für Mathematik, Technische Universität Dortmund, D-44227 Dortmund, Germany}

\begin{abstract}
Let $X \subset \mathbb{P}^{n}$ be an unramified real curve with $X(\mathbb{R}) \neq \emptyset$. If $n \geq 3$ is odd, Huisman \cite{huisnon} conjectured that $X$ is an $M$-curve and that every branch of $X(\mathbb{R})$ is a pseudo-line. If $n \geq 4$ is even, he conjectures that $X$ is a rational normal curve or a twisted form of a such. Recently, a family of unramified $M$-curves in $\mathbb{P}^3$ that serve as a counterexample to the first conjecture was constructed in \cite{maximallyw}. In this note we construct another family of counterexamples that are not even $M$-curves. We remark that the second conjecture follows for generic curves of odd degree from the de Jonquières formula.
\end{abstract}


\end{frontmatter}

\section{Introduction}
In this note, a real (algebraic) curve $X$ is assumed to be smooth, geometrically integral, embedded into the real projective space $\mathbb{P}^{n}$ for some $n \in \mathbb{N}^*$, and such that the set of real points $X(\mathbb{R})$ is nonempty. Since we assumed the set of real points $X(\mathbb{R})$ to be nonempty, it inherits the structure of an analytic manifold, which decomposes into a finite number of connected components, which are called $\emph{(real) branches}$ of $X$. We call such a branch an \emph{oval} if its homology class in $H_{1}(\mathbb{P}^{n}(\mathbb{R}), \mathbb{Z}/2)$ is trivial, and a \emph{pseudo-line} otherwise. By Harnack's Inequality \cite{harnack}, we know $s \leq g+1$, where $s$ is the number of branches and $g$ is the genus of $X$. We say that $X$ is an \emph{$M$-curve} if it has the maximum number of branches, i.e., if $s=g+1$. The curve $X$ is called $\emph{nondegenerate}$ if there is no real hyperplane $H \subset \mathbb{P}^{n}$ such that $X \subset H$. A nondegenerate curve $X$ is called $\emph{unramified}$ if, taken any real hyperplane $H$, we have
\[ \text{wt}(H \cdot X) \leq n-1,\]
whereby the weight of the intersection divisor $H\cdot X$ is defined to be
\[\text{wt}(H \cdot X) \coloneqq \text{deg}\left(H\cdot X - \left(H \cdot X\right)_{\text{red}}\right),\]
i.e., the degree of the difference between the latter and the reduced divisor (which contains each point of $H\cap X$ with multiplicity exactly one). Otherwise, it is called \emph{ramified}. For example, if $X\subset\mathbb{P}^{2}$ is a plane curve, two different types of ramification can occur: \emph{Bitangents}, lines that are tangent to the curve in two different points, and \emph{flex lines}, lines that intersect the curve in some point with multiplicity (at least) three. A result due to Klein \cite{kleinformel} implies that the numbers $f$ and $b$ of real flex lines and real bitangents respectively to a real plane curve of degree $d$ satisfy the inequality $$f+2b\geq d(d-2).$$In particular, every real plane curve of degree $d\geq3$ is ramified. The analogous statement fails in odd-dimensional projective spaces: Huisman \cite[Thm.~3.1]{huisnon} shows that in any $\mathbb{P}^{n}$ with $n\geq3$ odd there are unramified $M$-curves of arbitrary large genus and degree. However, in the same article he makes the following conjecture:
\begin{con}[Conjecture 3.4 in \cite{huisnon}]\label{con:odd}
Let $n \geq 3$ be an odd integer and $X \subset \mathbb{P}^{n}$ be an unramified real curve.
Then $X$ is an $M$-curve and each branch of $X$ is a pseudo-line, i.e., it realizes the nontrivial homology class in $H_{1}(\mathbb{P}^{n}(\mathbb{R}), \mathbb{Z}/2)$.
\end{con}
It was pointed out in \cite[Rem.~2.7]{rigid} that \Cref{con:odd} is false: The \emph{maximally writhed links} constructed in \cite{maximallyw} are instances of unramified $M$-curves in $\mathbb{P}^3$ but not all of them consist of pseudo-lines only.

We remark that \Cref{con:odd} has an interesting link to the following question on totally real divisor classes:
\begin{que}\label{qu:totreal}
Given a real curve $X$, determine the smallest natural number $N(X) \in \mathbb{N}^*$ such that any divisor of degree at least $N(X)$ is linearly equivalent to a totally real effective divisor, i.e., a divisor whose support consists of real points only.
\end{que}
It was shown by Scheiderer \cite[Cor.~2.10]{scheiderer} that such a number $N(X)$ exists. \Cref{qu:totreal} was studied by Huisman \cite{huison} and Monnier \cite{monnier}, but it seems challenging to obtain results for curves with few branches. However, assuming \Cref{con:odd} to be true, Monnier \cite[Thm.~3.7]{monnier} established a new bound for $N(X)$ at least for $(M-2)$-curves (i.e., $s=g-1$). We remark that Huisman \cite{huisun} has shown \Cref{con:odd} under more restrictive assumptions (namely nonspecial linearly normal curves having ``many branches and few ovals''). The main objective of this note is to construct families of counterexamples to \Cref{con:odd} living on the quadric hyperboloid in $\mathbb{P}^{3}$ that are, in contrast to those from \cite{maximallyw}, not even $M$-curves.

\bigskip
 \noindent \textbf{Acknowledgements.}
We would like to thank Daniel Plaumann for the very helpful discussions that initiated this project.

\section{Unramified curves in three-space}
This section includes two explicit counterexamples and a method for constructing families of counterexamples. First, we present a simple idea for generating an unramified rational curve of degree $4$ in the three-space. Since it has even degree, its single branch is an oval, and hence contradicts the second part of \Cref{con:odd} though being an $M$-curve (as stated in the first part).
\begin{ex}\label{exp:firstex}
The coefficients of the fourth power of a nonzero homogeneous real linear form in two variables define a rational normal curve
\[\varphi: \mathbb{P}^{1} \longrightarrow \mathbb{P}^{4}, \left[a:b \right] \longmapsto \left[a^4: 4a^3b: 6a^2b^2 : 4ab^3:b^4\right].\]
The real hyperplane $H\coloneqq \mathcal{V}_{+}(-2x_0+x_1-x_2+x_3-2x_4)$ does not intersect the real image points of $\varphi$; the intersection divisor
\[H \cdot \text{im}(\varphi) = 2\rho^{\sigma}\]
with $\rho=\varphi\left(\left[e^{i\frac{\pi}{3}}:1\right]^{\sigma}\right)$
consists of a complex conjugate pair with multiplicity two. The linear projection from the orthogonal complement of $H$ defines a map
\[\psi : \mathbb{P}^{1} \longrightarrow \mathbb{P}^{3}, [x:y] \longmapsto \left[x^4+2x^3y:x^4-2x^2y^2:x^4+2xy^3:-x^4+y^4\right].\]
The image is a real curve in $\mathbb{P}^{3}$ of genus $0$ and degree $4$, hence the only branch is not a pseudo-line. The curve is indeed unramified, since intersecting with an arbitrary real hyperplane is equivalent to finding the homogeneous roots of a polynomial 
\[ \lambda_0(x^4+2x^3y) + \lambda_1(x^4-2x^2y^2) + \lambda_2(x^4+2xy^3) + \lambda_3 (-x^4+y^4), \]
which cannot be a single real point with multiplicity $4$; thus the weight of the intersection with an arbitrary hyperplane is at most $2$. The curve is cut out by a quadric of the Segre type and three cubics.
\end{ex}

The following more general construction produces a series of unramified curves in $\mathbb{P}^{3}$ of any degree and any even genus. In contrast to \cite[Rem.~2.7]{rigid}, we stress that these curves are constructed to be far away from being $M$-curves; they consist of exactly one branch.

\begin{cons}\label{con:main} Let $p,q \in \mathbb{R}[t]$ be \emph{strictly interlacing} polynomials both of degree $d \in \mathbb{N}^*$, i.e., all complex zeros of $p$ and $q$ are real and between each two consecutive zeros of $p$ there is exactly one zero of $q$, see \cite[\S 6.3]{rahsch}. The graph $Y$ of the fraction $p(t)/q(t)$ has the following shape in the real plane:

\begin{figure}[H]
	\centering
  \includegraphics[width=0.4\textwidth]{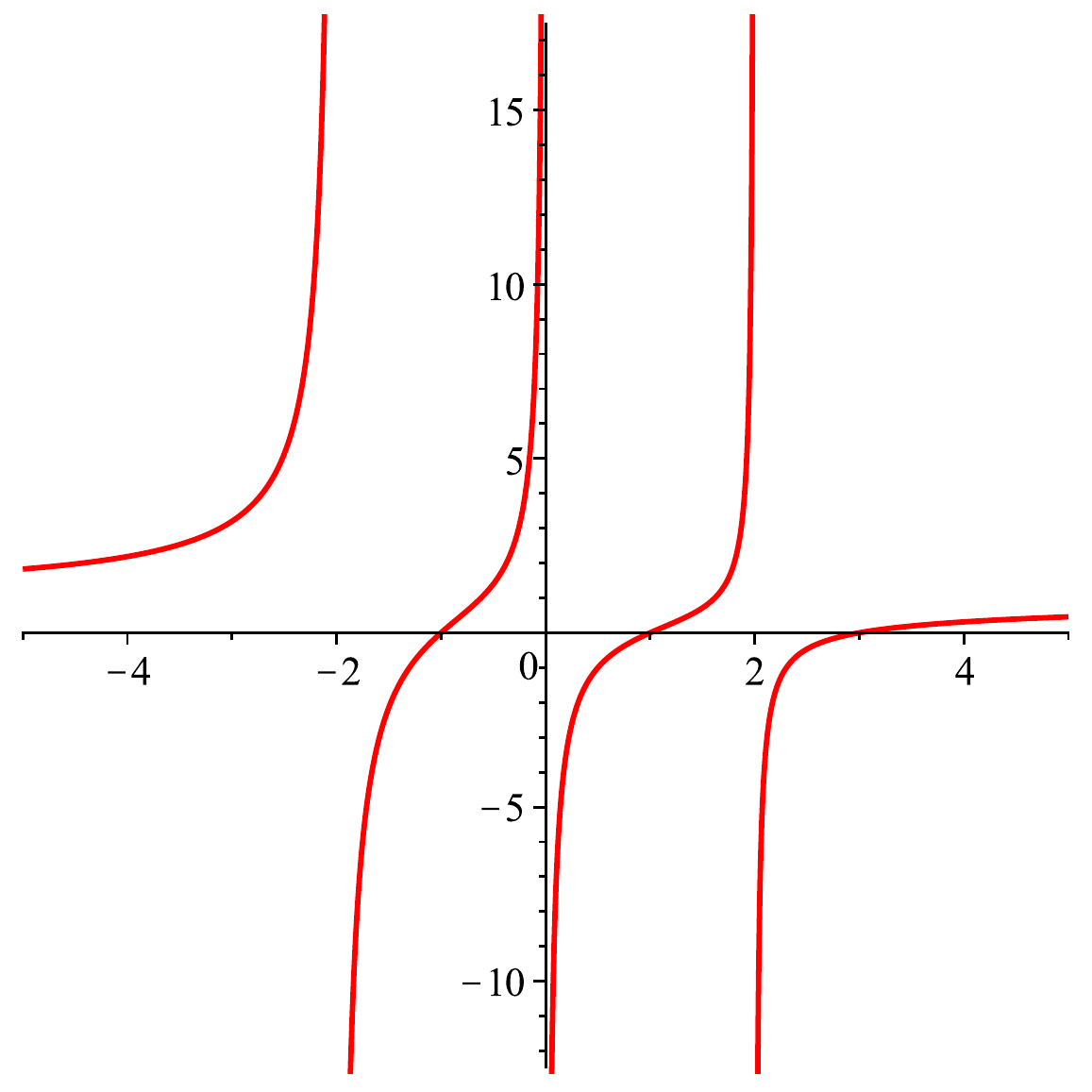}
	\caption{$p(t)=(t^2-1)(t-3)$ and $q(t)=(t^2-4)t$}
	\label{figure1}
\end{figure}

\noindent Let $X \coloneqq \overline{Y}$ denote the Zariski closure of that graph in $\mathbb{P}^{1} \times \mathbb{P}^{1}$, i.e., $X$ is the image of the map
\[\varphi: \mathbb{P}^{1} \longrightarrow \mathbb{P}^{1} \times \mathbb{P}^{1}, \left[x:y\right] \longmapsto \left(\left[x :y\right], \left[ Q(x,y) : P(x,y)\right]\right),\]
where $P$ and $Q$ denote the homogenizations of $p$ and $q$ (e.g., $P(1,t)=p(t)$). The curve $X$ realizes the homology class $\left(1,d\right)$ in $H_{1}\left(\Sigma, \mathbb{Z}\right)$ where $\Sigma=\left(\mathbb{P}^{1}\times \mathbb{P}^{1}\right)(\mathbb{R})$. Indeed, since it is the graph of a function, any vertical real line intersects $X$ in exactly one point and it follows for example from the Hermite--Kakeya Theorem \cite[Thm.~6.3.8]{rahsch} that every real horizontal line intersects $X$ in $d$ simple real points. After we embed $\mathbb{P}^{1} \times \mathbb{P}^{1}$ to $\mathbb{P}^{3}$ via the Segre embedding, the real intersection with a hyperplane has homology class $(1,1)$ (or $(1,-1)$) and thus intersects $X$ in 
\[\det\begin{pmatrix}
1 & 1\\
1 & d
\end{pmatrix}= d-1 \hspace{5px} \left(\text{or} \hspace{5px} \det\begin{pmatrix}
1 & -1\\
1 & d
\end{pmatrix}= d+1 \right)\]
different real points. In particular, the induced embedding of $X$ to $\mathbb{P}^{3}$ is unramified. Thus up to now, we have constructed a family of examples of unramified, rational curves of degrees $d+1$ in $\mathbb{P}^{3}$. Our curve from Example \ref{exp:firstex} also arises in that way for $d=3$ (up to a change of coordinates). In order to obtain curves of positive genus, let $G$ denote the Zariski closure of
\[\mathcal{V}\left(\prod_{j=1}^{e}(t^2+t_j)\right) \subset \mathbb{A}^{2} \subset \mathbb{P}^{1} \times \mathbb{P}^{1}\]for pairwise distinct positive $t_j$
embedded as in the previous case via 
\[\mathbb{A}^{2} \longrightarrow \mathbb{P}^{1} \times \mathbb{P}^{1}, (s,t) \longmapsto \left([1:s],[1:t]\right).\]
The union $Z \coloneqq X \cup G$ realizes the class $(2e+1,d)$ in the Picard group of $\mathbb{P}^{1} \times \mathbb{P}^{1}$, i.e., it is defined by a bihomogeneous polynomial $F$ of bidegree $(2e+1,d)$. We have $Z(\mathbb{R})=X(\mathbb{R})$. We now observe that every real hyperplane $H$ intersects $Z$ in $e$ different complex conjugate pairs, one for each pair of lines, as well as in at least $d-1$ (for $H$ of the type $(1,1)$) or in exactly $d+1$ (for $H$ of the type $(1,-1)$) real points. Thus $H$ intersects $Z$ in at least $d+2e-1$ different points. This remains true for any sufficiently small perturbation $F_{\epsilon}$ of $F$ with a smooth zero set $Z_{\epsilon}$ and we can conclude that
\[ \deg\left( H \cdot Z_{\epsilon} \right)_{\text{red}} \geq d+2e-1, \]
hence the weight of $H \cdot Z_{\epsilon}$ is at most two for all real hyperplanes $H$, i.e., $Z_{\epsilon}$ is unramified. We further note that in $\mathbb{P}^3$ the curve $Z_{\epsilon}$ has degree $d+2e+1$ and for $\epsilon$ sufficiently small it still realizes the homology class $(1, d)$ in $H_{1}(\Sigma, \mathbb{Z})$. Finally, the constructed real curve has the one branch coming from $X$, which is an oval or pseudo-line depending on the parity of $d$, and no other branches. The genus is $2e\cdot(d-1)$, see \cite[Ch. V, Exp.~1.5.2]{hartshorne}.
\end{cons}

\begin{rem}
By \cite[Thm.~1]{maximallyw}, the rational curves of degrees $d+1$ in $\mathbb{P}^3$ from above (i.e., before adding line-pairs) are maximally writhed knots. Thus for $e=0$ our counterexamples agree with those from \cite[Rem.~2.7]{rigid}. Letting $d=1$ and $e>1$ we obtain a family of unramified rational curves of even degree that are not maximally writhed.
\end{rem}

\begin{rem}
Each of our examples $Z_{\epsilon}$ from Construction \ref{con:main} intersects any real horizontal line only in real points. This implies that we obtain a morphism $f:Z_{\epsilon}\longrightarrow\mathbb{P}^1$ which is \emph{totally real} in the sense that $f^{-1}(\mathbb{P}^1(\mathbb{R}))=Z_{\epsilon}(\mathbb{R})$. This implies that $Z_{\epsilon}$ is of \emph{type I} in the sense that $Z_{\epsilon}(\mathbb{C})\setminus Z_{\epsilon}(\mathbb{R})$ is not connected. We don't know whether there are also counterexamples to \Cref{con:odd} that are not of type I.
\end{rem}

\begin{rem}
We remark that the counterexamples of positive genus do not stay unramified after changing the base field to $\mathbb{C}$. In fact, the only nondegenerate curves that are unramified over $\mathbb{C}$ are the rational normal curves. This follows from the Pl\"ucker formulas, see \cite[Exc.~I.C-14]{acgh}.
\end{rem}

\begin{ex}
We consider the case $d=2$ and $e=1$. In this case our resulting curve $Z_{\epsilon}$ has genus two and degree five. In particular, the curve is embedded to $\mathbb{P}^3$ via a complete and nonspecial linear system. In the light of \cite[Thm.~2]{huisun}, which says that \Cref{con:odd} holds for nonspecial rational and elliptic curves, this is a minimal counterexample with respect to the genus and being nonspecial.
\end{ex}

Finally, we compute explicit equations for the case $d=3$ and $e=1$.

\begin{ex}
For the following computations we used Macaulay2 \cite{M2}.
We take $p(t)=(t^2-1)(t-3)$ and $q(t)=(t^2-4)t$; the graph of the fraction is shown in \Cref{figure1}. The closure in $\mathbb{P}^{1} \times \mathbb{P}^{1}$ consists of points coming from $\mathbb{A}^{2}$, i.e., from $\left(\left[1:t\right],\left[1: p(t)/q(t)\right]\right) \hspace{3px} \text{for} \hspace{3px} t \in \mathbb{A}^{1}\setminus \lbrace -2,0,2 \rbrace ,$ together with four points at infinity connecting appropriate end points of the graph in $\mathbb{A}^{2}$. In the real projective three-space, the ideal of the resulting curve $X$ is defined by the Segre quadric and three cubics. The closure of the complex conjugate pair of lines (setting $e=1$ and $t_1=1$) in $\mathbb{P}^{3}$ is
\[G \coloneqq \mathcal{V}_{+}\left(\langle x_{0}x_3-x_1x_2, x_0^2+x_1^2,x_2^2+x_3^2 \rangle\right).\]
The intersection of the ideal of $X$ with the ideal of $G$ is the ideal generated by $q\coloneqq x_{0}x_3-x_1x_2$ and the following polynomial (denoted by $h$ henceforth):
\[3x_0^3+3x_0x_1^2-x_0^2x_2-3x_0x_2^2+x_2^3+4x_0^2x_3-x_0x_1x_3+4x_1^2x_3-x_2^2x_3-3x_0x_3^2+x_2x_3^2-x_3^3.\]
We deform $h$ by $p\coloneqq x_0^3+x_1^3+x_2^3-x_3^3$. The result $Z_{\epsilon}\coloneqq \mathcal{V}_{+}\left(q, h+ \epsilon p\right)$ for a small $\epsilon$ is an explicit example of an unramified curve of degree six, genus four, i.e., a \emph{canonical curve}, with exactly one oval and of type I. Note that an unramified canonical curve of genus four must be of type I with one oval because otherwise it would have a real tritangent plane \cite[Prop.~5.1]{grossharris}.
\end{ex}

\begin{figure}[H]
	\centering
   \includegraphics[width=0.4\textwidth]{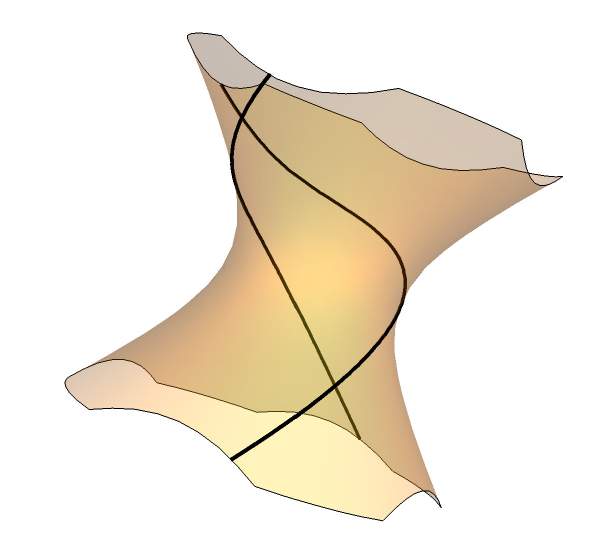}
 	\caption{The curve $Z_{10^{-5}}$ living on the hyperboloid}
 	\label{figure2}
\end{figure}

\section{Unramified curves in even-dimensional spaces}
For the sake of completeness, we remark that Huisman \cite{huisnon} has also studied embeddings of $M$-curves into even-dimensional projective spaces. For certain such embeddings he computes the exact number of real inflection points which turns out to be nonzero for \emph{positive genus}. Here by an inflection point of a curve $X\subset\mathbb{P}^n$ we mean a point $P\in X$ for which there is a hyperplane $H$ intersecting $X$ in $P$ with multiplicity at least $n+1$.
In particular, a curve having an inflection point is ramified. In his article, he conjectures the following general statement:

\begin{con}[Conjecture 4.6 in \cite{huisnon}]\label{con:even}
Let $n \geq 4$ be an even integer and $X \subset \mathbb{P}^{n}$ be an unramified real curve. Then $X$ is a rational normal curve or a twisted form of a rational normal curve (i.e., after changing the base field to $\mathbb{C}$, the curve $X_{\mathbb{C}}$ is a rational normal curve).
\end{con}
Again, we remark that Huisman \cite{huisun} has shown \Cref{con:even} under more restrictive assumptions (namely nonspecial linearly normal curves having ``many branches and few ovals'').

\begin{rem}
We would like to point out that this conjecture is in fact true for generic curves of \emph{odd degree}. Indeed, by the de Jonquières formula \cite[p.~359]{acgh} a generic nondegenerate curve in $\mathbb{P}^{2n}$ of degree $2d+1$ and genus $g$ has$$(2d+1)\cdot(2n+1)+2n\cdot(2n+1)\cdot (g-1)$$ complex inflection points. Since this is clearly an odd number, there must be a real inflection point. Therefore, the curve is ramified.
\end{rem}
\nocite{bcr}
\center{\bibliography{References1}}

\end{document}